\documentclass[debug]{rmaa}
\usepackage{paralist}

\usepackage{psfrag,color}

\usepackage[latin1]{inputenc}
\usepackage{amsmath}
\usepackage{url}
\usepackage{graphicx}
\usepackage{natbib}
\usepackage{float}
\title{Location and stability of the triangular points in the triaxial elliptic restricted three-body problem}
\author{M. Radwan \altaffilmark{1} and Nihad S. Abdel-Motelp \altaffilmark{2}}

\altaffiltext{1}{Astronomy; Space Science and Metreology Department, Faculty of Science, Cairo University, 12613 Giza, Egypt.}

\altaffiltext{2}{Astronomy Department, National Research Institute of Astronomy and Geophysics (NRIAG), 11421 Helwan, Cairo, Egypt.}
\shortauthor{Radwan \& Abdel-Motelp}
\shorttitle{Location and stability of $L_{4,5}$ in the triaxial ER3BP}
\fulladdresses{
	\item M. Radwan: Astronomy; Space Science and Metreology Department, Faculty of Science, Cairo University, 12613 Giza, Egypt (mradwan@sci.cu.edu.eg).
	\item Nihad S. Abd El Motelp: Astronomy; Astronomy Department, National Research Institute of Astronomy and Geophysics (NRIAG), 11421 Helwan, Cairo, Egypt.
	(nihad.saad@nriag.sci.eg).}
\listofauthors{M. Radwan, \& N . S. Abd El Motelp}
\indexauthor{Radwan, M.}
\indexauthor{Abd El Motelp, N. S}

\abstract
{The main goal of the present paper is to evaluate the perturbed locations and investigate the linear stability of the triangular points. We studied the problem in the elliptic restricted three body problem frame of work. The problem is generalized in the sense that the two primaries are considered as triaxial bodies. It is found that the locations of these points are affected by the triaxiality coefficients of the primaries and the eccentricity of orbits. Also it is observed that the stability regions depend on the involved perturbations. In addition to this we studied the periodic orbits in the vicinity of the triangular points.}

\resumen{}
\addkeyword{Celestial Mechanics}
\addkeyword{Gravitation}
\addkeyword{Methods: Analytical}
\addkeyword{Methods: Numerical}
\addkeyword{Planets and Satellites: Dynamical Evolution And Stability}
\begin{document}
	\maketitle
	\section{Introduction}
	\label{Int}
	The restricted three-body problem (R3BP) is the most significant problem in celestial mechanics. This due to its wide applications in space dynamics and solar system dynamics. The knowledge of the R3PB is crucial in almost all space applications. We can define this problem as a dynamical system in which an infinitesimal body moves in the gravitational field of two massive bodies (primaries). The two primaries move about their common centre of mass in either circular or elliptic orbits, which leads to the famous circular or elliptic restricted problems, respectively \citep{10.1088/0034-4885/77/6/065901}.
	In the elliptic restricted three-body problem ER3BP, it is assumed that the smaller primary moves around the more massive one in an elliptical orbit. The influence of the infinitesimal mass on the motion of the primaries is negligible. This significantly simplified model has more applications in space dynamics than the more general problem.
	\par{Despite the fact that the RTBP is not integrable, five specified solutions in the rotating frame are obtained. These solutions correspond to equilibrium positions. Three of these equilibrium points, known as collinear points, are located on the line joining the two massive primaries in the rotating reference frame. The remaining two points, called triangular points, form equilateral triangles with the massive primaries. The determination of the location of these libration points are of great interest to different space applications \citep{10.1016/j.newast.2018.10.005}. In general the equilateral equilibrium points are conditionally stable, while the collinear libration points are usually unstable \citep{10.1140/epjp/i2016-16365-s}.}
	\par{Most of the celestial bodies are of irregular nature and cannot be assumed spherically permanent in the RTBP, because their stability of movement is influenced by their shape. The planets of the solar system and their satellites are assumed as extended bodies and cannot be regarded as spheres. Furthermore, the effect of eccentricity on the orbits is significant. So that in order to achieve an acceptable accuracy, we have to take into account that the orbits of most celestial bodies are elliptic rather than circular \citep{10.1140/epjp/i2016-16365-s}.}
	\par{From the theory of small oscillations around the triangular points two types of periodic orbits can be obtained; short-period orbits with periods approximately the same as the orbital period of the primaries and long-period orbits with periods many time that of the primaries. The importance of investigating these periodic orbits is situated in their significant appearing in nature, they provide us with important information about orbital resonances, spin orbits and also to approximate quasi-periodic trajectories (\citet{10.1007/s11038-013-9415-5}, \citet{10.1007/s10509-009-0038-2}, \citet{10.17485/IJST/v13i32.601}, \citet{10.2478/amns.2020.2.00022}).}
	\par{A number of studies have been performed on the equilateral triangular points of the restricted three-body problem. \citet{10.1007/s10509-010-0545-1} investigated the existence of periodic orbits about the triangular points, namely $L_4$ and $L_5$, in the RTBP when the primaries are considered as triaxial and oblate spheroid bodies, besides perturbations due to Coriolis and centrifugal forces. The authors found that long and short periodic orbits exist around these points and their eccentricities and periods are influenced by the involved perturbations. 
		\citet{10.1088} studied the stability of triangular points in the ER3BP, when both oblate primaries emit light energy. The locations of $L_4$ and $L_5$ are shifted away from the line joining the two primaries than in the classical case. The stability region decreases and increases with variability in eccentricity, oblateness and radiation pressures.
		\citet{10.1007/s10509-014-1818-x} investigated the stability of infinitesimal motions around the equilateral points. They considered the more massive primary as a radiant source and the less massive one as a triaxial body.
		\citet{10.12988/nade.2016.6863}, \citet{10.2298/SAJ1795047Z} studied the locations and linear stability of triangular points under the effects of oblateness and triaxiality of the primaries, plus small relativistic perturbations. They found that the locations of the triangular points are affected by the considered perturbations. Numerical explorations are carried out to show the effect of perturbations on the position and stability.
		\citet{10.1140/epjp/i2016-16365-s} investigated the positions and stability of the libration points when both primaries are taken as oblate spheroids with oblateness up to the fourth zonal harmonic $J_4$. They found that both harmonic coefficients $J_2$ and $J_4$, eccentricity and semi-major axis have destabilizing tendencies and as a result the size of the region of stability decreases with an increase in the parameters involved.}
	\par{Recently, \citet{1937-1632_2019_4-5_703} studied analytically the existence and the linear stability of the libration points assuming the primaries as triaxial bodies and Euler angles are equal to specified values. The authors proved that the positions and the stability of the triangular points change according to the effect of the triaxiality of the primaries. Furthermore, they presented the solution of long and short periodic orbits for stable motion.}
	\citet{10.1134/S156035472002001X} revisited the circular spatial restricted problem. The authors presented some new results in the light of the concept of Lie stability. This paper studied the positions and stability of the well-known triangular points, in the elliptic restricted problem frame of work. The triaxiality of the primary bodies are taken into consideration. In addition, we studied the periodic orbits around the triangular points. To examine the present problem we carried out several numerical explorations.
	%%%%%%%%%%%%%%%%%%%%%%%%%%%%%%%%%%%%%%%%%%%%%%%%%%%%%%%%%
	\section{Dynamical model}
	Let $m_1$, $m_2$ and $m_3$ be the masses of the primary, secondary, and the third infinitesimal body, respectively. The third body moves under the effect of the triaxial primaries, but it does not affect their motion. The two primaries move around their common centre of mass in elliptic orbits. We consider a rotating coordinate system $(x,y,z)$ with origin at the centre of mass of the primaries. The primaries are permanently located on the x-axis, which is the line joining them. When investigating the restricted three-body problem it is better to choose a system of normalized units. So we assume that the distance between the primaries is unity and we choose also the unit of time to make the gravitational constant unity. The total mass of the primaries are considered as $1$, i.e. $m_1+m_2=1$. In this normalized system we adopt $\mu=\dfrac{m_2}{m_1+m_2}$ as a dimensionless mass parameter. According to the above consideration, the equations of motion of the third infinitesimal body moves under the effect of the triaxial primaries located at $(\mu, 0)$, $(1-\mu, 0)$ are given as \citet{10.2478/amns.2020.2.00022}.
	%%%%%%%%%%%%%%%%%%%%%%%%%%%%%%%%%%%%%%%%%%%%%%%%%%%%%%%%%
	\begin{equation}
		\begin{aligned}
			\ddot{x}-2\:n\:\dot{y} &=\frac{\partial U}{\partial x}, \qquad
			\ddot{y}+2\:n\:\dot{x} &=\frac{\partial U}{\partial y},
		\end{aligned}
	\end{equation}
	with
	%%%%%%%%%%%%%%%%%%%%%%%%%%%%%%%%%%%%%%%%%%%%%%%%%%%%%%%%%%%%%%
	\begin{equation}
		\begin{split} 
			U & =\dfrac{n^2}{2}\:\left[\left( 1-\mu\right)\:r_1^2+\mu\:r_2^2\right]+
			\dfrac{\left(1-\mu\right)}{r_1}+\dfrac{\mu}{r_2}+\\
			&\dfrac{\left( 1-\mu\right)\:\left(2\:\sigma_1- \sigma_2\right)}{2\:r_1^3}
			-\dfrac{3\:\left( 1-\mu\right)\:\left(\sigma_1-\sigma_2\right)\: y^2}{2\: r_1^5}\\
			&+\dfrac{\mu\:\left(2\:\gamma_1-\gamma_2\right)}{2\:r_2^3}-
			\dfrac{3\:\mu\:\left(\gamma_1-\gamma_2\right)\:y^2}{2\:r_2^5}.
		\end{split}
		\label{sultan2018}
	\end{equation}
	%%%%%%%%%%%%%%%%%%%%%%%%%%%%%%%%%%%%%%%%%%%%%%%%%%%%%%%%%%%%%%%%%%%%%
	\begin{equation}
		\begin{split}
			r_1 & = \sqrt{\left( x+\mu\right)^2+y^2},\\
			r_2 & = \sqrt{\left( x+\mu-1\right)^2+y^2},\\
		\end{split} 
		\label{distance}
	\end{equation}
	\begin{equation}
		\begin{split}
			n & = \frac{1}{a \left(1-e^2\right)} \biggl(1+\frac{3}{2}\:
			\biggl(\left(2\:\sigma_{1}-\sigma_{2}\right)\\
			&+\left(2\:\gamma_{1}-\gamma_{2}\right)\biggr)\biggr).
			\label{meananomaly}
		\end{split}
	\end{equation}
	Where $U$ is the restricted three-body potential; 
	$n$ is the mean motion of the primaries; $r_1$ and $r_2$ are the distances of the infinitesimal mass from the primaries; 
	$\sigma_{1}$, $\sigma_{2}$ and $\gamma_{1}$, $\gamma_{2}$ are the triaxiality coefficients of the bigger and smaller primaries, respectively.
	%%%%%%%%%%%%%%%%%%%%%%%%%%%%%%%%%%%%%%%%%%%%%%%%%%%%%%%%%%%%%%%%%%%%%%
	\section{Locations of triangular points}
	
	The libration points represent stationary solutions of the restricted three-body problem. The positions of these points can be obtained by setting all components of relative velocity and acceleration equal to zero, i.e. we found them by setting 
	$\dot{x}=\dot{y}=0=\ddot{x}=\ddot{y}$. Consequently, the locations of these points can be found by solving simultaneously the nonlinear equations $U_x=U_y=0$, and 
	$y \neq0$. Setting
	\begin{align*}
		A_{\gamma} = & \gamma_1-\gamma_2,\qquad A_{\sigma} = \sigma_1-\sigma_2,\\
		A_{\alpha} = & 2\:\gamma_1-\gamma_2,\quad A_{\beta} = 2\:\sigma_1-\sigma_2,
	\end{align*}
	we have
	\begin{equation}
		\begin{split}\label{ux}
			U_x & = -\mu\:\left(-1+x+\mu\right)\biggl[ \frac{1}{r_2^3}+\frac{3}{2\:r_2^5}\:A_\alpha-\frac{15\:y^2}{2\:r_2^7}\:A_\gamma\biggr
			]\\
			&+\frac{x}{a\:\left(1-e^2\right)}+
			\left(1-\mu\right)\:\left(x+\mu\right)\biggl[\frac{1}
			{r_1^3}-\frac{3}{2\:r_1^5}\:A_\beta\\
			&+\frac{3}{2\:a\:\left(1-e^2\right)}\:
			\left(A_\alpha+A_\beta\right)-
			\frac{15\:y^2}{2\:r_1^7}\:A_{\sigma}\biggr]
		\end{split}
	\end{equation}
	%%%%%%%%%%%%%%%%%%%%%%%%%%%%%%%%%%%%%%%%%%%
	\begin{equation}
		\begin{aligned}
			U_y & = -y\:\biggl[\frac{\mu}{r_2^3}+\frac{1-\mu}{r_1^3}+
			\frac{3\:\mu}{2\:r_2^5}\:A_{\alpha}+
			\frac{3\:\left(1-\mu\right)}{2\:r_1^5}\\
			&\:A_{\beta} -\frac{1}{a\:\left(1-e^2\right)}\:\left[1+\frac{3}{2}\left(A_\alpha+A_\beta\right)\right]-\\
			&\frac{15\:y^2}{2}\:\left[\frac{\mu\:A_\gamma}{r_2^7}+
			\frac{\left(1-\mu\right)}{r_{1}^{7}}\:A_\sigma\right]-\\
			& \frac{3\:\mu}{r_2^5}\:A_\gamma+\frac{3\:\left(1-\mu\right)}{r_1^5}\:A_\sigma\biggr],\\
		\end{aligned}
		\label{uy}
	\end{equation}
	%%%%%%%%%%%%%%%%%%%%%%%%%%%%%%%%%%%%%%%%%%%
	Ignoring the triaxiality perturbations due to the primary bodies yields the equilateral solution of the classical restricted three-body problem i.e. $r_1=r_2=1$, then it may be reasonable to consider that the locations of these points are the same as given by classical problem but perturbed by small terms, i.e.   
	\begin{equation}
		r_i = 1+\delta_i,\qquad\delta_i<<1,\qquad \left(i=1,2\right).
		\label{perturbeddisctance}
	\end{equation}
	Substituting this assumption \eqref{perturbeddisctance} into \eqref{distance} then solving for $x$ and $y$ and retaining terms up the first order in the small quantities $\delta_i$ we get:
	%%%%%%%%%%%%%%%%%%%%%%%%%%%%%%%%%%%%%%
	\begin{equation}
		\begin{aligned}
			x & = \frac{1}{2}\left(2\:\delta_1-2\:\delta_2-2\:\mu +1\right),\\
			y & =\pm\frac{1}{2} \sqrt{\left(3+4\:\left(\delta_1+\delta_2\right)\right)}.\\
			\label{lcoordinate}
		\end{aligned}
	\end{equation}
	Substituting the values of $r_1$, $r_2$, $x$ and $y$ into equations \eqref{ux} and \eqref{uy} and ignoring the terms higher than order one in $\delta_i$. Thus we obtain the following two simultaneous equations in $\delta_1$ and $\delta_2$. 
	%%%%%%%%%%%%%%%%%%%%%%%%%%%%%%%%%%%%%%%%%
	\begin{equation}
		A_1+A_2\:\delta_1+A_3\:\delta_2=0,\qquad B_1+B_2\:\delta_1+B_3\:\delta_2=0.
	\end{equation}
	The corresponding solution is:
	\begin{equation}
		\delta_1 = -\frac{A_3\:B_1-A_1\:B_3}{A_3\:B_2-A_2\:B_3},\qquad  {A_3\:B_2-A_2\:B_3}\neq0,
	\end{equation}
	\begin{equation}
		\begin{split} 
			\delta_2 = -\frac{-A_2\:B_1+A_1\:B_2}{A_3\:B_2-A_2\:B_3},\qquad {A_3\:B_2-A_2\:B_3}\neq0,
		\end{split}
	\end{equation}
	%%%%%%%%%%%%%%%%%%%%%%%%%%%%%%%%%%%%%%%%%%%%%%%%%%%%
	where:
	\begin{equation}
		\begin{split}
			A_1 & = -\frac{1}{2}+\mu+\frac{3}{4}\:\left(\mu\: A_\alpha-A_\beta+\mu\: A_\beta\right)+\\
			&\frac{3}{4\:a\: \left(1-e^2\right)}\:\left(A_\alpha+A_\beta\right)+ 
			\frac{3}{2\:a \left(1-e^2\right)}\:\\
			&\left(1-\mu\:\left(2+3\:A_\alpha+3\:A_\beta\right)\right)+\frac{45}{16}\:\\
			&\left(-\mu\:\left(A_\gamma+ A_\sigma\right)
			+A_\sigma\right),
		\end{split}
	\end{equation}
	%%%%%%%%%%%%%%%%%%%%%%%%%%%%%%%%%%%%%%%%%%%%%%%%%%%%%%%%%%%%%%%
	\begin{equation}
		\begin{split}
			A_2 & =\frac{1}{2}-\frac{1}{a \left(1-e^2\right)}
			\left(1+\frac{3}{2}\left(A_\alpha+A_\beta\right)\right)+
			\frac{3}{2}\\
			&\biggl(1-\mu\:A_\alpha+
			\frac{3}{2}\:\left(1-\mu\right)\: A_\beta-\frac{55}{8}\:\left(1-\mu\right)\:A_\sigma\\
			&+\frac{5}{4}\:\mu\:A_\gamma\biggr),
		\end{split}
	\end{equation}
	%%%%%%%%%%%%%%%%%%%%%%%%%%%%%%%%%%%%%%%%%%%%%%5
	\begin{equation}
		\begin{split}
			A_3 & = 1+\frac{3}{2}\:A_{\beta}+\frac{1}{a\: \left(1-e^2\right)}\:\left(1+\frac{3}{2}\:\left(A_\alpha+A_\beta\right)\right)-\\
			&\frac{15}{8}\:A_\sigma-
			\frac{3}{2}\:\mu\left(1+\frac{3}{2}A_\alpha+A_\beta-
			\frac{55}{8}-\frac{5}{4}\:A_{\sigma}\right),
		\end{split}
	\end{equation}
	%%%%%%%%%%%%%%%%%%%%%%%%%%%%%%%%%%%%%%%%%%%%%%%%%%
	\begin{equation}
		\begin{split}
			B_1 & = -\frac{
				\sqrt{3}}{2}\:\biggl(1-\frac{1}{a\:\left(1-e^2 \right)} \biggl(1+\frac{3}{2}\:\biggl(A_\alpha+\\
			&A_\beta\biggr)\biggr)
			-\frac{3}{2}\:\biggl(A_\beta-\frac{7}{4}\: A_{\sigma}+
			\mu\:\biggl(A_\alpha-A_\beta-\\
			&\frac{7}{4} A_\gamma+
			\frac{7}{4}\:A_{\sigma}\biggr)\biggr)\biggr),
		\end{split}
	\end{equation}
	%%%%%%%%%%%%%%%%%%%%%%%%%%%%%%%%%%%%%%%%%%%%%%%%%%%%%%
	\begin{equation}
		\begin{split}
			B_2 & = \frac{7}{2\:\sqrt{3}}+\frac{1}{\sqrt{3}\:a\:\left(1-e^2 \right)}\:\biggl(1+
			\frac{3}{2}\biggl(A_\alpha+\\
			&A_\beta\biggr)\biggr)+
			\frac{\sqrt{3}}{2}\:\biggl(\frac{13}{2}\:A_\beta-\frac{121}{8}\:A_\sigma
			-\\
			&\mu\:\biggl(3+A_\alpha+\frac{13}{2}\:A_\beta-
			\frac{37}{4}\:A_\gamma-\frac{121}{8}A_\sigma \biggr)
			\biggr),
		\end{split}
	\end{equation}
	%%%%%%%%%%%%%%%%%%%%%%%%%%%%%%%%%%%%%%%%%%%%%%%%%
	\begin{equation}
		\begin{split}
			B_3 & = -\frac{1}{\sqrt{3}}\left(1+\frac{1}{a
				\left(-1+e^2\right)}\right)+\\
			&\frac{3\sqrt{3}}{2}
			\biggl(\frac{37}{4}\:A_{\sigma}-
			\frac{A_\beta}{3}+
			\mu\: \biggl(1+\frac{13}{2}\: A_\alpha+\\
			&A_\beta-\frac{121}{8} \: A_\gamma-\frac{37}{4}\:A_\sigma\biggr)-
			\frac{\left(A_\alpha+
				A_\beta\right)}{a\:\left(-1+e^2\right)}\biggr).
		\end{split}
	\end{equation}
	%%%%%%%%%%%%%%%%%%%%%%%%%%%%%%%%%%%%%%%%%%%%%%%%%%
	\section{Graphical representations (location)}                      
\begin{figure}[!t]
		\includegraphics[width=\columnwidth]{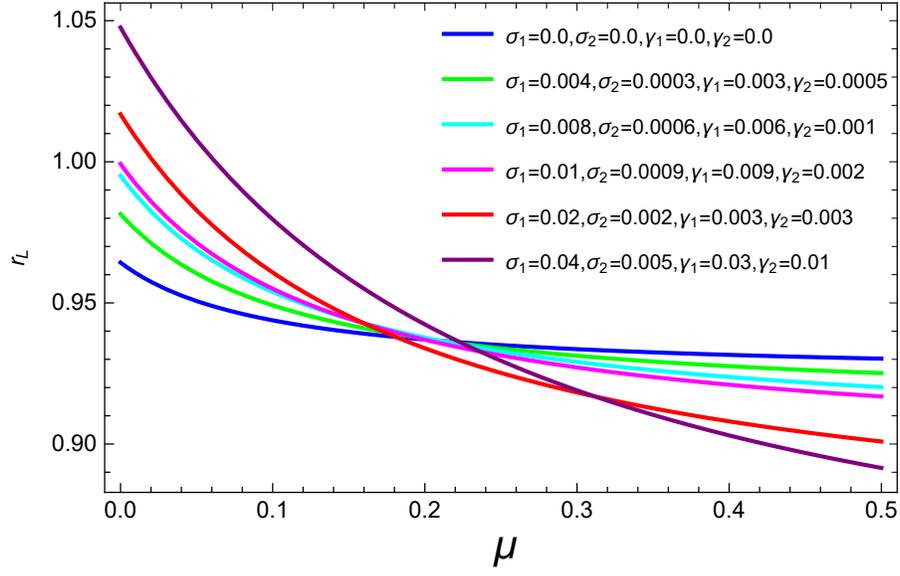}
		\caption{The locations of $L_{4,5}$ with different triaxiality coefficients and eccentricity $e=0.2$}
		\label{fig:triaxialparameter}
	\end{figure}

\begin{figure}[!t]
		\includegraphics[width=\columnwidth]{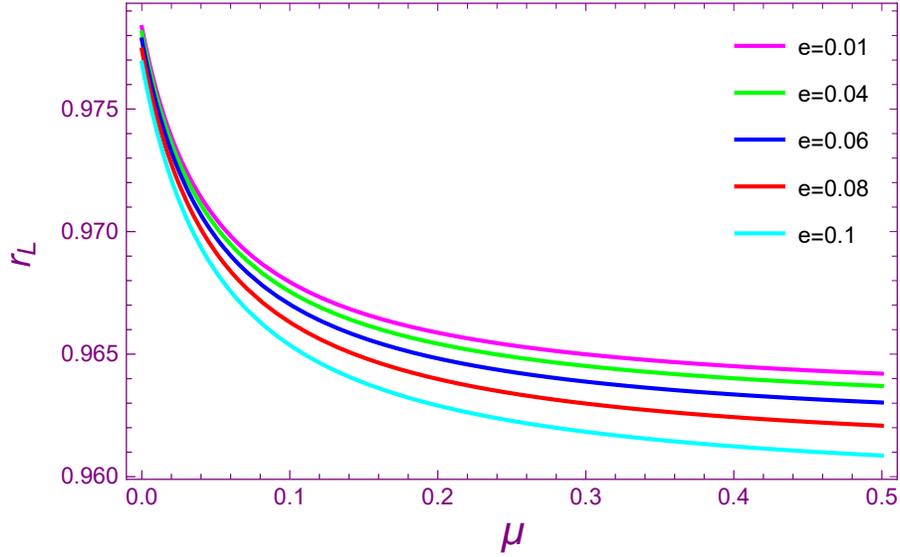}
		\caption{The locations of $L_{4,5}$ with different eccentricities and triaxiality  coefficients $\sigma_1=0.0002$, $\sigma_2=0.0005$, $\gamma_1=0.003$ and $\gamma_2=0.002$.}
		\label{fig:changingineccentricity}
	\end{figure}

	In Fig.\ref{fig:triaxialparameter}-\ref{fig:changingineccentricity} we plotted the locations $r_{L_{4,5}}$ of triangular points in the elliptic restricted three body problem, taking into consideration triaxiality and eccentricity effects. Fig.\ref{fig:triaxialparameter} represents the variation of $r_{4,5}$ against the mass ratio $\mu$. The calculations are carried out for several dynamical models with different coefficients $\gamma_i$ and $\sigma_i$; and constant eccentricity of primaries orbit $e=0.2$. As is clear from the figure, the locations begin to decrease sharply until we reached the value $\mu \approx 0.15$, then it decrease slowly. We also observed that the variation in the locations $ r_{L_{4,5}}$ is small, because the triaxiality coefficients are very small. Fig.\ref{fig:changingineccentricity} depicts the variation of $r_{L_{4,5}}$ against the mass ratio $\mu$ with constant triaxiality coefficients. We applied our case study for a range of eccentricities, the eccentricity $e$ takes the values $0.01, 0.04, 0.06, 0.08$ and $0.1$. 
	We noticed that the change in curves is nearly the same for all considered cases, this because the chosen values of $e$ are small. The dynamics is approximately the same but with different sizes of perturbations. We expect that the higher the eccentricity the larger the perturbation in the positions of triangular points.
	%%%%%%%%%%%%%%%%%%%%%%%%%%%%%%%%%%%%%%%%%%%%%%%%%%%%%%%%%%%%%%%%
	%%%%%%%%%%%%%%%%%%%%%%%%%%%%%%%%%%%%%%%%%%%%%%%%%%%%%5
	%%%%%%%%%%%%%%%%%%%%%%%%%%%%%%%%%%%%%%%%
	\section{Stability of Triangular points}
	The dynamical models which describe the restricted problem are very difficult. The great complexity of these models lead us to focus our attention on the linearized dynamics, because in this case we obtain simplified mathematical expressions that could be handled easily. The location of the third body is displaced a little from the equilibrium point due to the considered perturbations. If the resultant motion of this body is a rapid departure from the proximity of the point, we can call such a position of equilibrium point an unstable one, if the body just oscillates about the point, it is said to be a stable position \citep{10.1016/j.newast.2018.10.005}. To examine the possible motion of the infinitesimal body in the proximity of triangular points, the equations of motion should be linearized.
	Let the position of the triangular points be ($x_{\circ}$,$y_{\circ}$) and let the third body be displaced to the point ($x_{\circ}+\xi$, $y_{\circ}+\eta$), where $\xi$ and $\eta$ are small displacement in ($x_{\circ}, y_{\circ}$), so that the linearized equations can be written as:
	%%%%%%%%%%%%%%%%%%%%%%%%%%%%%%%%%%%%%%%%%
	\begin{equation}
		\begin{split}
			\ddot{\xi}-2\:n\:\dot{\eta} & = \frac{1}{a\left(1-e^2\right)}\biggl(U_{xx}^{L_{4,5}}\:\xi+
			U_{xy}^{L_{4,5}}\: \eta\biggr)\\
			\ddot{\eta}-2\:n\:\dot{\xi} &= \frac{1}{a\left(1-e^2\right)}
			\biggl( U_{yy}^{L_{4,5}}\:\eta+U_{xy}^{L_{4,5}}\:\xi\biggr)
			\label{linerizedequations}
		\end{split}
	\end{equation}
	Where $U_{xx}^{L_{4,5}}$ denotes the second derivative of $U$ with respect to $x$ computed at the stationary solution and the rest of the derivatives can be defined similarly. The characteristic equation corresponding to Equ.\eqref{linerizedequations} may be written as
	\begin{equation}
		\begin{aligned}
			\lambda^4+\left(4\:n^2-U_{xx}^{L_{4,5}}-
			U_{yy}^{L_{4,5}}\right)\: \lambda^2+U_{xx}^{L_{4,5}}\:U_{yy}^{L_{4,5}}-\left(U_{xy}^{L_{4,5}}\right)^2=0,
			\label{lamda4}
		\end{aligned}
	\end{equation}
	which may be written in the form
	\begin{equation}
		\lambda^4+U_{L_{i}}\: \lambda^2+\nu_{L_i}=0,
		\label{lamda}
	\end{equation}
	with
	\begin{equation}
		\begin{split}
			U_{L_{i}} & = \left(4\:n^2-U_{xx}^{L_{4,5}}-
			U_{yy}^{L_{4,5}}\right),\\
			\nu_{L_i} & =U_{xx}^{L_{4,5}}\: U_{yy}^{L_{4,5}}-\left(U_{xy}^{L_{4,5}}\right)^2,\\
			\label{lamda1}
		\end{split}
	\end{equation}
	the roots of equation \eqref{lamda1} are 
	\begin{equation}
		\lambda_{1,2}=\pm\frac{\sqrt{-U_{L_{i}}-
				\sqrt{U_{L_{i}}^2-4\:\nu_{L_{i}}}}}{\sqrt{2}},
		\label{lamda12}
	\end{equation}
	and
	\begin{equation}
		\lambda_{3,4}=\pm\frac{\sqrt{-U_{L_{i}}+
				\sqrt{U_{L_{i}}^2-4\:\nu_{L_i}}}}{\sqrt{2}}.
		\label{lamda34}
	\end{equation}
	The possible solutions for equation \eqref{lamda4} will determine the stability of the equilibrium points.
	We have three possibilities: asymptotically stable when all the solutions are real and negative, unstable when at least one of the solutions is positive and real and the last one stable case if all solutions are pure imaginary \citep{9780126806502}. we will focus our attention only on the case which have an oscillatory stable solutions about triangular points.
	To investigate the stability/instability we first compute the partial derivatives required for equations 
	%%%%%%%%%%%%%%%%%%%%%%%%%%%%%%%%%%%%%%%
	\eqref{lamda12} and \eqref{lamda34} as
	\begin{equation}
		\begin{split}
			U_{xx} & =\frac{1}{a\:\left(1-e^2\right)}\left(1+\frac{3}{2}
			\left(A_\alpha+A_\beta\right)\right)+\\
			&\mu\:
			\biggl(\frac{3\:\left(\mu+x-1\right)^2}{r_2^5}-
			\frac{1}{r_2^3}+
			\frac{3}{2}\:A_{\alpha} 
			\biggl(\frac{5\:\left(\mu+x-1\right)^2}{r_2^7}\\
			&-\frac{1}{r_2^5}\biggr)-
			\frac{15}{2} y^2 A_{\gamma}
			\left(\frac{7}{r_2^9}\:\left(\mu+ x-1\right)^2-
			\frac{1}{r_2^7}\right)\biggr)\\
			&+\left(1-\mu\right)\:\left(\frac{3}{2}\:A_{\beta}\:\biggl(\frac{5\left(\mu +x\right)^2}{r_1^7}-
			\frac{1}{r_1^5}\right)+\\
			&\left(\frac{3 \left(\mu +x\right)^2}{r_1^5}-
			\frac{1}{r_1^3}\right)
			-\frac{15}{2} y^2 A_{\sigma}
			\biggl(\frac{7\left(\mu+x\right)^2}{r_1^9}
			-\frac{1}{r_1^7}\biggr)\biggr),
		\end{split}
	\end{equation}
	%%%%%%%%%%%%%%%%%%%%%%%%%%%%%%%%%%%%%%%%%%%%%%%%%%%%%%%
	\begin{equation} 
		\begin{aligned}
			U_{yy} & = \frac{1}{a\:\left(1-e^2\right)}\:\left(1+
			\frac{3}{2}\left(A_{\alpha}+A_{\beta}\right)\right)\\
			&+\mu\:\left(3\:A_\gamma\:\left(\frac{10\:y^2}{r_2^7}-\frac{5\: y^2}{2}\left(\frac{7\: y^2}{r_2^9}-\frac{1}{r_2^7}\right)-\frac{1}{r_2^5}\right)\right)\\
			&+\left(1-\mu\right)\:\biggl(\frac{3\: y^2}{r_1^5}-
			\frac{1}{r_1^3}+\frac{3}{2}\:\left(\frac{5\:y^2}{r_1^7}-
			\frac{1}{r_1^5}\right)\:A_\beta\\
			&-3\:A_\sigma\:\left(\frac{1}{r_1^5}+5\:y^2\left(\frac{-2}{r_1^7}+\frac{7\:y^2}{2\: r_1^9}-\frac{1}{2\:r_1^7}\right)\right)\\
			&+\frac{3\:A_\alpha}{2}\:
			\left(\frac{5\: y^2}{r_2^7}-\frac{1}{r_2^5}\right)+
			\frac{3\:y^2}{r_2^5}-\frac{1}{r_2^3}\biggr),
		\end{aligned}
	\end{equation}
	%%%%%%%%%%%%%%%%%%%%%%%%%%%%%%%%%%%%%%%%%%%%%%%
	\begin{equation}
		\begin{aligned}
			U_{xy} & = y\:\biggl(\mu\:
			\left(\mu +x-1\right)\:\biggl(-\frac{105\: y^2 A_{\gamma}}{2\: r_2^9}+\frac{15\:A_{\alpha}}{2\:r_2^7}\\
			&+\frac{15\:A_{\gamma }}{r_2^7}+
			\frac{3}{r_2^5}\biggr)+
			y\: \left(1-\mu\right)\:\left(\mu+x\right)\:
			\biggl(\frac{3}{r_1^5}\\
			&-\frac{105\:y^2\:A_\sigma}{2\:r_1^9}+
			\frac{15\:A_{\beta}}{2\:r_1^7}+
			\frac{15\:A_\sigma}{r_1^7}\biggr)\biggr).
		\end{aligned}
	\end{equation}
	%%%%%%%%%%%%%%%%%%%%%%%%%%%%%%%%%%%%%%%%%%%%%%%%%%%%%%%%%%%%%%%%%%%%%%%%%
	\section{Graphical representations (stability)}
\begin{figure}[tb]
		\includegraphics[width=\columnwidth]{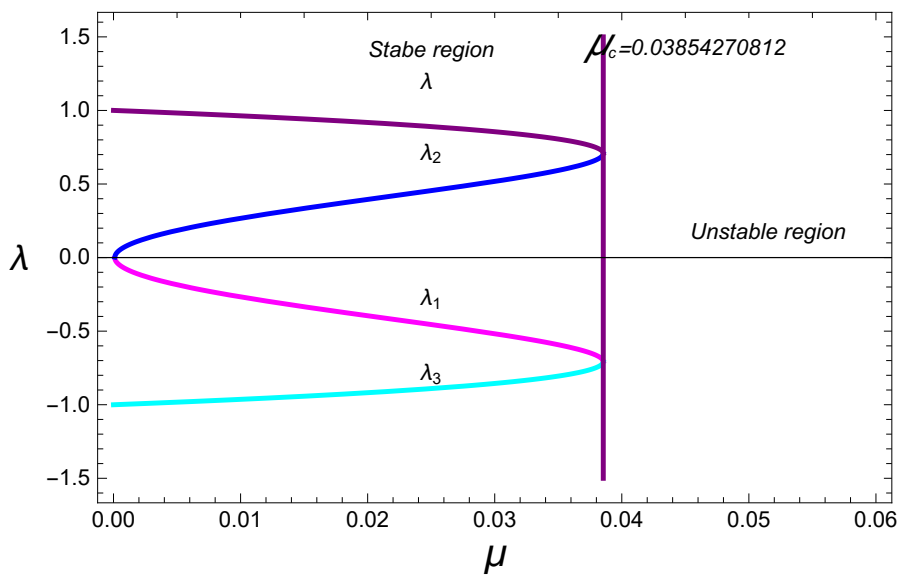}
		\caption{The imaginary part of the eigenvalues versus mass parameter in the classical case, $\sigma_i = \gamma_i=0$.}
		\label{fig:stablec}
	\end{figure}

\begin{figure}[!t]
		\includegraphics[width=\columnwidth]{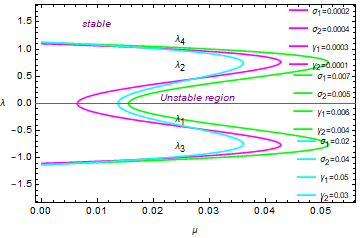}
		\caption{The imaginary part of the eigenvalues versus mass parameter with eccentricity $e=0.01$ and different triaxial parameters.}
		\label{fig:striaxial}
	\end{figure}

\begin{figure}[!t]
		\includegraphics[width=\columnwidth]{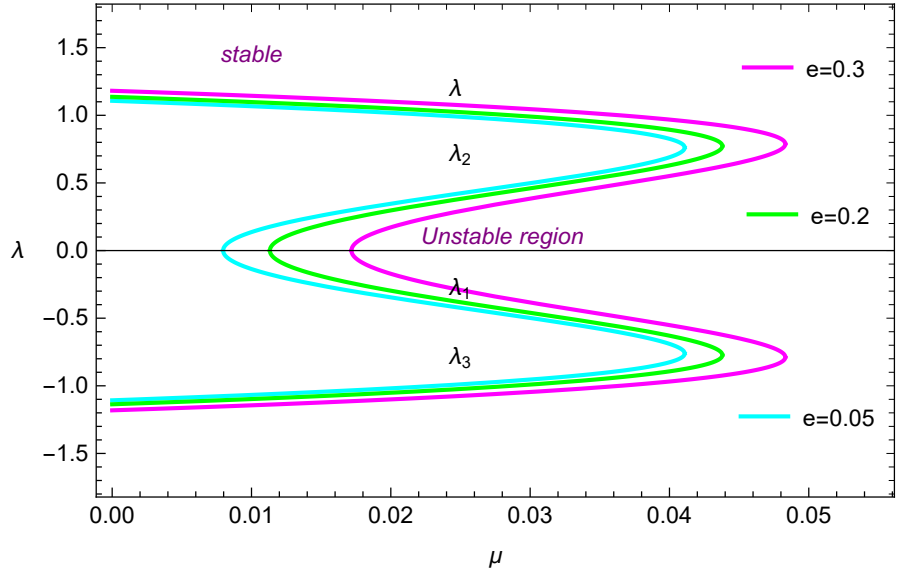}
		\caption{The imaginary part of the eigenvalues versus mass parameter 
			with different eccentricities and triaxial parameters $\sigma_1=0.001$, $\sigma_2=0.004$, $\gamma_1=0.003$ and $\gamma_2=0.002$.}
		\label{fig:seccentricity}
	\end{figure}

	Figs.~\ref{fig:stablec}-~\ref{fig:seccentricity} illustrate the stability regions taking into consideration the triaxiality coefficients and the eccentricity of the orbits. Fig.~\ref{fig:stablec} represents the classical case i.e. when $\sigma_{i}=0$, $\gamma_{i}=0$ and $e=0$ as well. 
	Fig.~\ref{fig:striaxial} depicts the stability regions for dynamical systems with different triaxiality coefficients, while Fig.~\ref{fig:seccentricity} shows these regions for different eccentricities and constant coefficients. We can see from Fig.~\ref{fig:striaxial} that each stability region is characterized by different values of the triaxiality coefficients. The size of the regions is not much different from that of the classical case, this because the perturbation due to the triaxiality coefficients is small. In Fig.~\ref{fig:seccentricity} the effect of different eccentricities on the size of the regions is obvious. Comparing with previous works e.g. \citet{10.1007/s10509-015-2308-5}, we found that the size of the regions is slightly different from that of the quoted paper due to the considered perturbations in each case. We expect that as the eccentricity increases the stability regions will be destroyed.
	Table~\ref{tbl:1} gives the values of critical mass for some selected dynamical models. These values are computed under effects of the included perturbations. We notice from the table that when the triaxiality coefficients ($\sigma_i$ and $\gamma_i$) are equal to zero, $\mu_c$ reduces to the critical mass value of the classical restricted problem (see Fig.~\ref{fig:stablec}). We also notice from Figs.~\ref{fig:striaxial},\ref{fig:seccentricity} that the triaxiality coefficients of the primaries change the range of stability.
	\par In addition, we can also investigate the behavior of the system by plotting the critical mass ratio against the eccentricity for different values of the triaxial parameters $\sigma_i$ and $\gamma_i$. We plotted  different combinations of triaxial parameters  and  a stable region has been observed. It is observed from Fig.6 that the increment in the eccentricity reduce the region of stability of the infinitesimal body about the triangular points. Figs.7 and 8 represent the variation of critical mass with the triaxial parameters when the eccentricity is constant. We can see from the two Figs. that increasing the triaxial parameters results in a destabilizing effect in the system.
	%-----------------------------------
	%-----------------------------------
	%%%%%%%%%%%%%%%%%%%%%%%%%%%%%%%%%%%%%%%%%%%%%%%%%%%%%%%%%%%%%%
	\begin{figure}[!t]
		\includegraphics[width=\columnwidth]{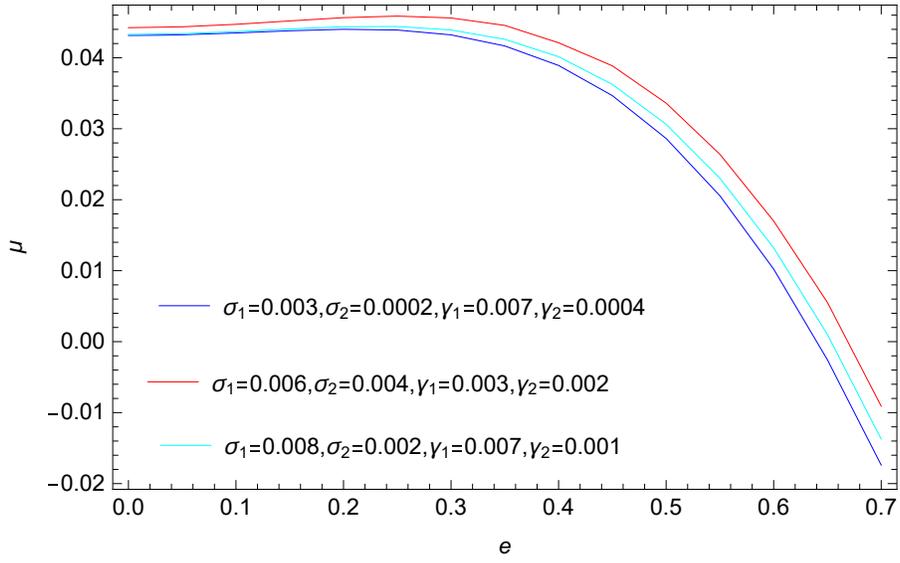}
		\caption{The critical mass ratio versus  eccentricity for different values of the triaxial parameters}
		\label{fig:eccentricityandmucritical}
	\end{figure}
	%%%%%%%%%%%%%%%%%%%%%%%%%%%%%%%%%%%%%%%%%%%%%
	\begin{figure}[!t]
		\includegraphics[width=\columnwidth]{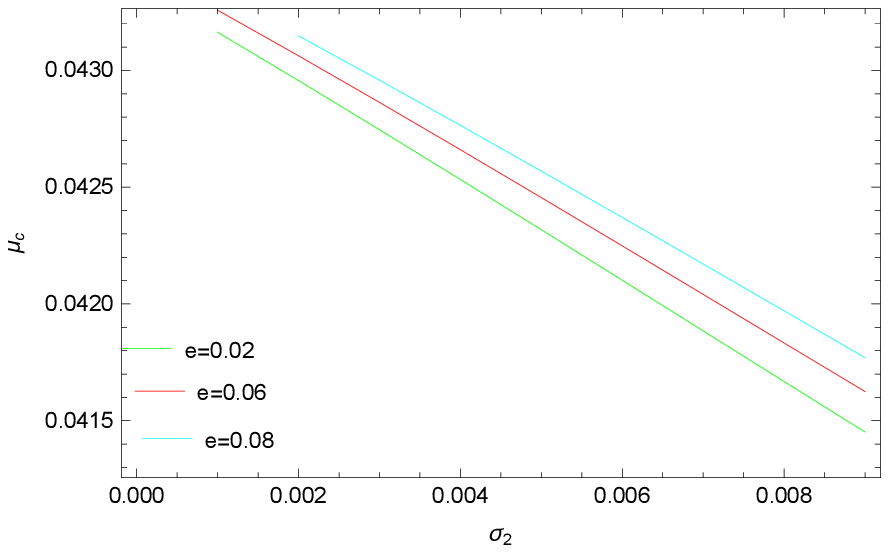}
		\caption{The critical mass ratio versus  triaxial parameter $\sigma_2$, for different values of eccentricity }
		\label{fig:sigma2_triaxiality_effecton_criticalmass}
	\end{figure}
	%%%%%%%%%%%%%%%%%%%%%%%%%%%%%%%%%%%%%%%%%%%%%%
	\begin{figure}[!t]
		\includegraphics[width=\columnwidth]{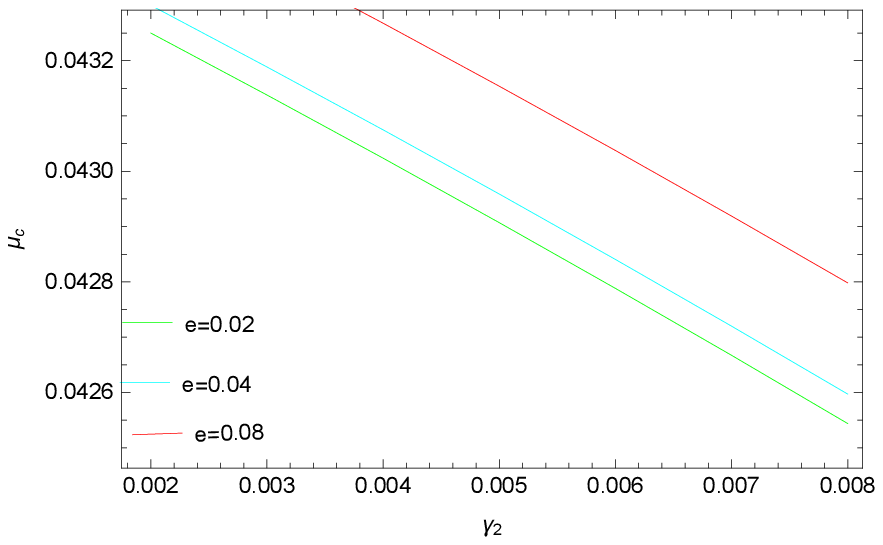}
		\caption{The critical mass ratio versus  triaxial parameter $\gamma_2$, for different values of eccentricity}
		\label{fig:gamma_2_triaxiality_effecton_criticalmass}
	\end{figure}
	%%%%%%%%%%%%%%%%%%%%%%
	\begin{table}[!t]\centering
		\small
		\caption{The values of the critical mass ratio for different dynamical models }
		\label{tbl:1}
		\tablecols{7}
		\label{table:1}
		\resizebox{\textwidth}{!}{%
			\begin{tabular}{llllllll}
				&$\sigma_1$ & $\sigma_2$ & $\gamma_1$& $\gamma_2$ & $e$ & $a$  & $\mu_c$\\
				\hline
				&  & &  &  & &  &  \\
				& 0 & 0  & 0 & 0 & 0& 1 & 0.0385209 \\
				&  & &  &  & &  &  \\
				&0.006 & 0.0004 & 0.03 & 0.002 & 0.2 & 0.92 & 0.0443621\\
				&0.006 & 0.0004 & 0.03 & 0.002 & 0.4 & 0.92 & 0.0401412\\
				&0.006 & 0.0004 & 0.03 & 0.002 & 0.6 & 0.92 & 0.0132093\\
				& & & & & & &\\
				& 0.003& 0.0002& 0.007& 0.0004& 0.02& 0.92& 0.0456447 \\
				& 0.003& 0.0002& 0.007& 0.0004& 0.04& 0.92& 0.0424109 \\
				& 0.003& 0.0002& 0.007& 0.0004& 0.07& 0.92& 0.00909276 \\
				&  & &  &  & &  &  \\
				& 0.009 & 0.002 & 0.008 & 0.002 & 0.02 & 0.92 & 0.042956\\
				& 0.009 & 0.004 & 0.008 & 0.002  & 0.02 & 0.92 & 0.0425331\\
				& 0.009 & 0.008 & 0.008 & 0.002 & 0.02 & 0.92 & 0.0416678\\			
				&  & &  &  & &  &  \\
				& 0.006 & 0.001 & 0.008 & 0.001 & 0.08 & 0.92 & 0.0466464\\
				& 0.006 & 0.001 & 0.008 & 0.005  & 0.08 & 0.92 & 0.0431538\\
				& 0.006 & 0.001 & 0.008 & 0.008 & 0.08 & 0.92 & 0.0427982\\			
		\end{tabular}}
	\end{table}
	%%%%%%%%%%%%%%%%%%%%%%%%%%%%%%%%%%%%%%%%%
	\section{Motion in proximity of triangular points}
	It is known that, periodic orbits represent the backbone in the study of dynamical systems of celestial mechanics. Since, in the range 
	$0 \leq \mu \leq \mu_c$, the roots of the characteristic polynomial of the present system are purely imaginary, thus we conclude that the motion around triangular points $L_4$ (or $L_5$) is stable and composed of two harmonic motions.
	The solutions for the perturbed motion about the triangular point $L_4$ (in the synodic reference frame) may be written as
	\begin{equation}
		X[t]  = \sum_{j=1}^{4}\:\alpha_j\:e^{\lambda_{j}\:t},
	\end{equation}
	\begin{equation}
		Y[t]  = \sum_{j=1}^{4}\:\beta_j\:e^{\lambda_{j}\:t},
	\end{equation}
	where the terms with coefficients $\alpha_j$ and $\beta_j$ represent the short and long periodic motion. The constants $\alpha_j$ and $\beta_j$ are not independent. Applying the initial conditions to the solutions of the system, then these solutions may be written in the form,
	\begin{equation}
		\begin{split}
			X[t] & =5.80204*10^{-5}\: \sin{\left(0.241582\: t\right)}-\\
			& 1.31958*10^{-5}\: \sin {\left(1.06221 \:t\right)}+\\
			& 3.19255*10^{-5}\: \cos {\left(0.241582 \:t\right)}-\\
			& 2.19255*10^{-5} \:\cos {\left(1.06221 \:t\right)}
		\end{split}
		\label{x1t}
	\end{equation}
	%==========================================
	\begin{equation}
		\begin{split}
			Y[t] & = 8.06288`*10^{-5}\:\cos{\left(0.241582`\: t\right)} - \\
			& 7.062876`*10^{-5} \:\cos{\left(1.06221`\: t\right)} - \\
			& 3.054769`*10^{-5}\:\sin{\left(0.241582`\: t\right)} + \\
			& 6.9475842 *10^{-5}\:\sin{\left(1.06221\: t\right)}\\
		\end{split}
		\label{y1t}
	\end{equation}
	%%%%%%%%%%%%%%%%%%%%%%%%%%%%%%%%%%%%%%%%%
	\begin{figure}[!t]
		\includegraphics[width=\columnwidth]{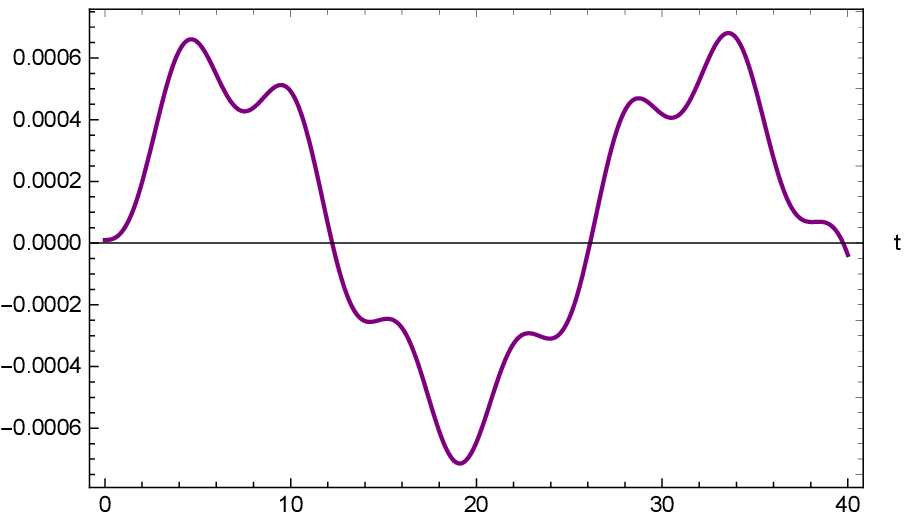}
		\caption{$x[t]$ as a function of time with triaxiality  coefficients 
			$\sigma _1 = 0.007$, $\sigma _2 = 0.002$, $\gamma _1 = 0.005$, 
			$\gamma _2 = 0.002$, $\mu=0.01$, $e = 0.01$ and $a=0.92$}.
		\label{fig:xt}
	\end{figure}
	%%%%%%%%%%%%%%%%%%%%%%%%%%%%%%%%%%%%%%%%%%%
	\begin{figure}[!t]
		\includegraphics[width=\columnwidth]{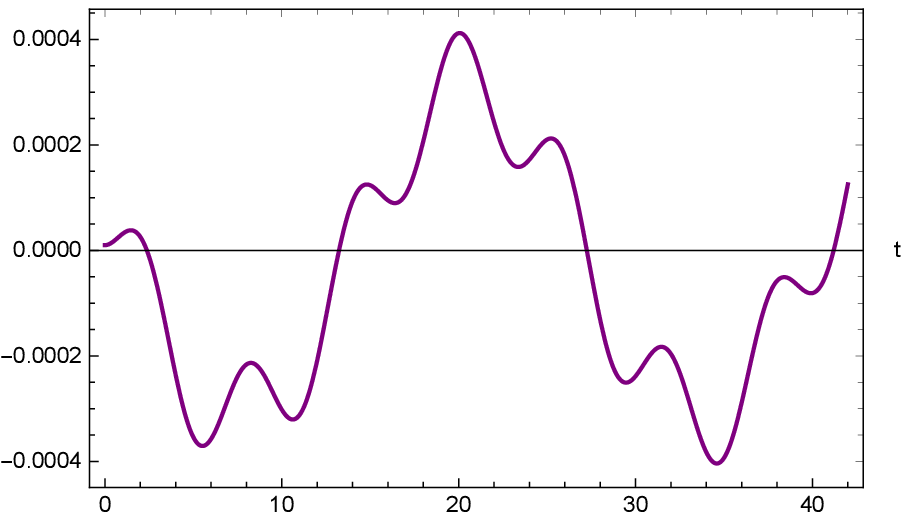}
		\caption{$y[t]$ as a function of time with triaxiality  coefficients
			$\sigma _1=0.007$, $\sigma _2=0.002$,$\gamma_1=0.005$, $\gamma_2=0.002$,  
			$\mu=0.01$, $e=0.01$ and $a = 0.92$.}
		\label{fig:yt}
	\end{figure}
	%=================================
	We may note, from equations \eqref{x1t} and \eqref{y1t}, that the solution is of the oscillatory type with fundamental periods of $\frac{2\:\pi}{0.241582`}$ and $\frac{2\:\pi}{1.06221}$ orbital periods of the orbiting body. Hence the body will remain in proximity of the triangular point $L_4$ and the motion is stable. These two different types of frequencies in the solution exist due to the fact that the resulting motion of the small body is composed of two types of motion. Firstly, a short-period motion with period close to the orbital period of the less massive primary and secondly a superimposed long-period one known as a libration around the point $L_4$. The two solutions are depicted in Figs.~\ref{fig:xt}, ~\ref{fig:yt}. The graph of the trajectory in the vicinity of equilibrium point $L_{4}$ are shown in Figs.~\ref{fig:comparisonl} and ~\ref{fig:comparisons} with magenta color. The graphs contain the short and the long-period motions, simultaneously. Fig.~\ref{fig:comparisonl} shows a comparison between the long-period motions in both the classical case (blue curve) and the present dynamical system (magenta curve), while Fig.~\ref{fig:comparisons} shows a comparison between the two systems in the short-period case. It is observed that, the trajectories (Magenta curve) in both Figs. are shifted from that of the classical case, this is due to the included perturbations.
	%========================================
	\begin{figure}[!t]
		\includegraphics[width=\columnwidth]{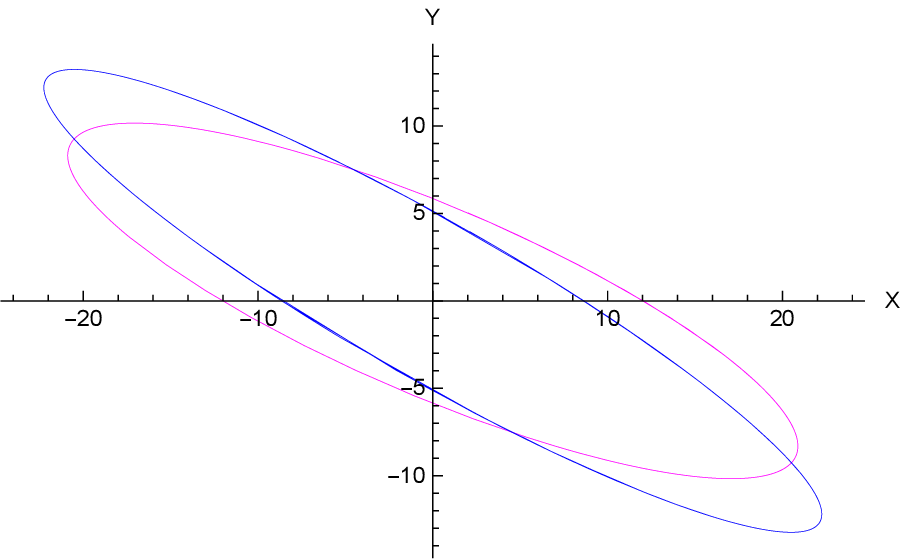}
		\caption{comparison between the long period motion in both classical case (blue curve) and the present dynamical system (Magenta curve).}
		\label{fig:comparisonl}
	\end{figure}
	%======================================
	\begin{figure}[!t]
		\includegraphics[width=\columnwidth]{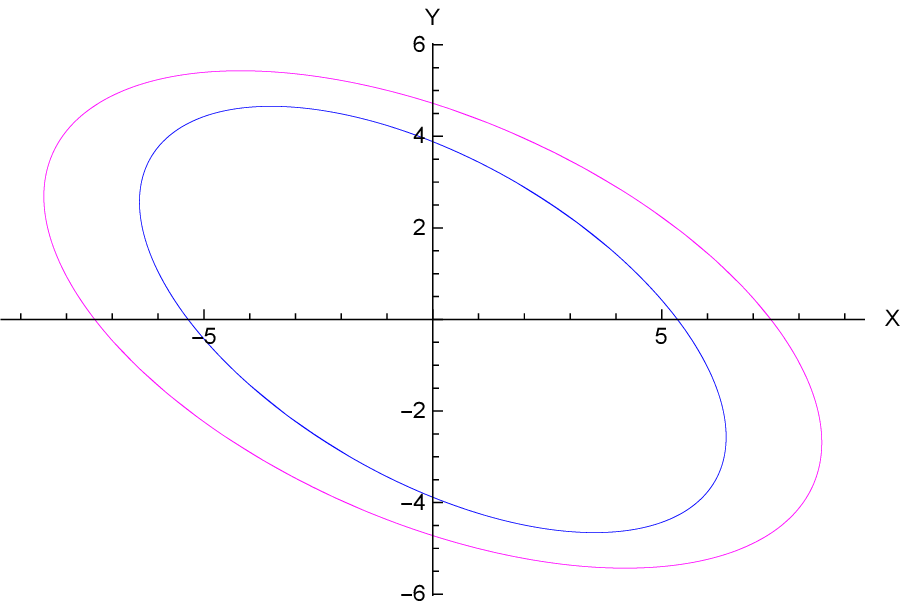}
		\caption{comparison between the short period motion in both classical case (blue) and the present dynamical system (Magenta).}
		\label{fig:comparisons}
	\end{figure}
	%%%%%%%%%%%%%%%%%%%%%%%%%%%%%%%%%%%%%%%%%
	
	%%%%%%%%%%%%%%%%%%%%%%%%%%%%%%%%%%%%%%%%%%%%%%%%%%%%%%%%%%%%
	\section{Conclusions}
	We have studied the existence and linear stability of the triangular points in the elliptic restricted three-body problem when the primaries are triaxial bodies. It is found that the locations and linear stability of the triangular points are affected by the triaxiality  coefficients of the primaries and the eccentricity of the orbit.
	It is also noticed that the size of the stability regions depend on the triaxiality  coefficients and the eccentricity of the orbits.
	In addition, we also studied the periodic orbits in the vicinity of the triangular points. It is seen that these orbits are elliptical. Further, it is observed that the shape of periodic orbits changed due to the triaxiality of the primaries and showed a deviation from the classical case.
	%References
	%% Please cite all reference entries in the article text using \cite or
	\clearpage
	%\bibliographystyle{plain}
	
	%\bibliography{refrence32}
	
\end{document}